\documentclass[12pt]{article}
\usepackage{amsmath}
\usepackage{graphicx,psfrag,epsf}
\usepackage{enumerate}
\usepackage{natbib}
\usepackage{url} 

\newcommand{\blind}{0}

\usepackage{amssymb}
\usepackage{tabularx,ragged2e,booktabs,caption}
\long\def\quote#1\par{{\bigskip\narrower{\openup -4 pt\small\noindent\noindent #1\par\bigskip}}\noindent}
\begin{document}

\def\spacingset#1{\renewcommand{\baselinestretch}%
{#1}\small\normalsize} \spacingset{1}


\if0\blind
{
  \title   
    {\LARGE\bf The likelihood principle does not entail a `sure thing', `evil demon' or `determinist' hypothesis}
    
  \author{Michael J. Lew
    \hspace{.2cm}\\
    Department of Pharmacology and Therapeutics\\University of Melbourne\\}
   
  \maketitle
} \fi

\if1\blind
{
  \bigskip
  \bigskip
  \bigskip
  \begin{center}
    {\LARGE\bf The likelihood principle does not entail a `sure thing', `evil demon' or `determinist' hypothesis}
    
\end{center}
  \medskip
} \fi

\bigskip
\begin{abstract}
The likelihood principle makes strong claims about the nature of statistical evidence but is controversial. Its claims are undermined by the existence of several examples that are assumed to show that it allows, with unity probability, domination of all other hypotheses by the uninteresting, determinist hypothesis that whatever happened had to happen. Such examples are generally assumed to be important obstacles to the application of the likelihood principle: they are counter-examples to the principle. A re-analysis of Birnbaum's 1969 ``counter-example'', demonstrates that the standardly reported analyses of such examples involves an inappropriate treatment of a nuisance parameter and that, when the nuisance parameter is adequately considered, there is no conflict between the evidential consequences of the likelihood principle and the intuitive evidential account of the problem. It also shows that the conclusion that the likelihood principle allows the determinist hypothesis to dominate with unity probability requires a misconception about the scope of the likelihood principle or an inappropriately specified statistical model. Whatever happened did \textit{not} have to happen.
\end{abstract}

\bigskip
\noindent%
\vfill

\newpage
\spacingset{1.45} 

\section{Introduction}
\label{sec:intro}
The likelihood principle offers an approach to utilizing the information in the whole of the likelihood function for inference and it appears to provide a principled basis for the assessment of statistical evidence. However, the likelihood principle is controversial.  It might seem to be a natural ingredient of most Bayesian approaches, and indeed it is a natural consequence of Bayes theorem and the product rule of probabilities \citep[pp.\ 250--251]{Jaynes:2003}, but Bayesians would usually not be comfortable with basing inference on likelihoods rather than probabilities. Of course, a Bayesian justification is not very compelling to frequentists who are effectively required to assume the likelihood principle is false because of its conflict with an evidential basis for frequentist approaches. Thus the likelihood principle has few friends. 

The likelihood principle has been proved to be entailed by the widely accepted principles of sufficiency and conditionality \citep{Birnbaum:1962}, but that proof has been challenged \citep{Mayo:2014} by a disproof that has itself been challenged in a different proof \citep{Gandenberger:2014}.  The existence, or not, of a proof may be a side-issue because it is possible that likelihood is a primitive, axiomatic postulate \citep{Jeffreys:1938cl, Fishe:1938, Edwards:1992, birnbaum:1969}. As Forster \& Sober pointed out ``If likelihood is epistemologically fundamental, then we should not be surprised to find that it cannot be justified in terms of anything more fundamental'' and that we could therefore expect only that the likelihood principle ``coincides with and systematizes our intuitions about examples''  \citep[p.\ 155]{ForsterSober:2004}. If agreement with our intuitions about examples is all we can expect then examples where the likelihood principle conflicts with our intuitions would stand as potentially fatal counter-examples.

The counter-example that is the subject of this paper is due to Alan \citet{birnbaum:1969}, who published the original proof of the likelihood principle in \citeyear{Birnbaum:1962}. He must have considered this example to be both important and persuasive, because he cites it in a letter to \textit{Nature} written as a response to a paper by A.W.F. Edwards promoting the use of likelihood-based scientific inferences \citep{Edwards:1969:nature}:

\quote{It is regrettable that Edwards's interesting article, supporting the likelihood and prior likelihood concepts, did not point out the specific criticisms of likelihood (and Bayesian) concepts that seem to dissuade most theoretical and applied statisticians from adopting them. \hfill ---\citep[p.\ 1033]{Birnbaum:1970}\\
}

\noindent Birnbaum also notes in that letter that he changed his mind about the likelihood between 1962 and 1969 and regretted not having made that clear at an earlier date. Edwards replied to Birnbaum's letter:

\quote{Birnbaum has retracted his former support for the likelihood approach to scientific inference after considering an example in which the approach would, it seems, always lead to the wrong conclusion. \hfill ---\citep{Edwards:1970:reply}
}

Birnbaum was not the only one-time proponent of the likelihood principle to have a change of heart as a consequence of apparent counter-examples. Ian Hacking, who first specified the law of likelihood \citep{hacking:1965} also wrote of such a counter-example as being important in the document where he notes that he had changed his mind concerning likelihood-based inference. The document was an extensive review of Edwards's 1972 book \textit{Likelihood}, and Hacking notes a variant of Birnbaum's example and writes that ``even on those increasingly rare days when I will rank hypotheses in order of their likelihoods, I cannot take the actual log-likelihood number as an objective measure of anything'' \citep[p.\ 137]{hacking:1972}. 


 Interestingly, the issue of \textit{The British Journal of Philosophy of Science} containing Hacking's review of Edwards's book also contains a (belated) review by G.A. Barnard of Hacking's own book of 1965, \textit{Logic of Statistical Inference}. Barnard was another of the leading proponent of likelihood-based inference, and his review notes the same problem as that raised by Birnbaum's example in the most general form:

\quote{The difficulty with Hacking's account is that it leads him to say (p. 89): `An hypothesis should be rejected if and only if there is some rival hypothesis much better supported than it is.' But there always is such a rival hypothesis, viz. that things just had to turn out the way they actually did.  \hfill ---\citet[p.\ 129]{Barnard:1972}

}

Given the seminal roles of Birnbaum, Hacking, Barnard and Edwards in the development of our understanding of the likelihood principle and the law of likelihood, their concerns about such counter-examples should be important to anyone with an interest in likelihood and statistical evidence.

%
%
%

\subsection{Structure}

The following sections of this paper consist of: (i) an introduction to likelihood functions and the likelihood principle, along with an example that contrasts probability and likelihood; (ii) a pair of simple examples which are similar to Birnbaum's problem but do not appear to throw up any conflict between intuition and the consequences of the likelihood principle; (iii) Birnbaum's example presented verbatim; (iv) a critique of Birnbaum's analysis and interpretation with two specific objections that invalidate the idea that Birnbaum's example is a counter-example to the likelihood principle; and (v) a discussion of the `whatever happened had to happen' hypothesis which demonstrates that that determinist hypothesis is not supported by a valid application of the likelihood principle.

\section{Likelihood and the likelihood principle}

Likelihood functions and the likelihood principle are rarely included in basic statistical textbooks and there is little consistency in form or naming of statements of the likelihood principle in the literature. Thus a brief discussion of likelihood and a clear statement of the likelihood principle is necessary to get this paper going. This section should be read by even experienced users of likelihood functions because it emphasizes some aspects of the likelihood principle that are usually tacit at best.

\subsection{Likelihood}
\label{L}

Assume a statistical model that has a parameter, $\theta$, that can take any value in the set $\Theta$. Such a model can be thought of as a meta-model consisting of a family of statistical models each of which assigns a single value to $\theta$. Each member of that family provides a probability distribution that gives the probability of any value, $x$, among the possible values, $X$, given the particular value of $\theta$. The probability distributions have values of $x$ as the horizontal axis, the `x-axis', and integrate to unity (or sum to unity where $X$ is discontinuous). Those probability distributions are the conventional probability distributions described and displayed by most introductory statistics textbooks. They are not likelihood functions.

Now, assume also that we wish to make inferences about values of the meta-model parameter $\theta$ and that we have observed a value, $x_{obs}$. The value $x_{obs}$ is special because it was actually observed, and thus the only relevant value of $x \in X$ is $x=x_{obs}$. The individual probability distributions of the family member models contain very little relevant information because $x_{obs}$ is only a single point, but if we plot the probabilities associated with $x_{obs}$ from all of the family member models together we get a function with an x-axis consisting of values of $\theta$. Such a function treats the observation as fixed and the parameter $\theta$ as variable, in contrast to the probability distributions described in the previous paragraph where $\theta$ is fixed and $x$ is the variable. The integral (or sum) of such a function can have any positive real value, which contrasts with a probability distribution where the integral (or sum) has to be unity. Such a function is a likelihood function.

Given the intimate relationship between probability distributions and likelihood functions it might be assumed that they are much the same thing. But that would be wrong because likelihoods are only defined (definable) to an arbitrary constant of proportionality, and that arbitrary scaling means that ``likelihood does not obey probability laws'' \citep[p.\ 17]{pawitan:2001:book}.

R.A. Fisher was the first to define likelihood:

\quote{
\textit{Likelihood}.---The likelihood that any parameter (or set of parameters)\footnote{Fisher uses the phrase ``set of parameters'' to mean a vector of parameters rather than a range of values for one parameter, as he makes clear where he deals with problems of estimation in that same paper.} should have any assigned value (or set of values) is proportional to the probability that if this were so, the totality of the observations should be that observed. \hfill --- \citep[p.\ 310]{Fisher:1922}
}

Thus we can write:

\begin{align}
\label{LfnDef}
& L(\theta \,|\, x_{obs}) \propto P(x_{obs} \,|\, \theta), \; \text{for all }\theta \text{ in } \Theta\\
& L(\theta \,|\, x_{obs}) = c.P(x_{obs} \,|\, \theta), \; \text{for all }\theta \text{ in } \Theta
\end{align}
\noindent where $c$ is an arbitrary constant.  Usually the ``given $x_{obs}$'' part is assumed and so the likelihood $L(\theta \,|\, x_{obs})$ is written as just $L(\theta)$. It is not unusual to see equations like \ref{LfnDef} but with `$=$' in place of the `$\propto$' without the provision of an arbitrary constant, and that is unfortunate because likelihoods and probabilities are not equivalent.

\subsection{The likelihood principle}
\label{LP}

In the paper where he sets out his alleged counter-example, Birnbaum describes two parts of the likelihood principle:

\quote{
Another general concept of statistical evidence, the \textit{likelihood concept} (often called the \textit{likelihood principle}), consists of two parts. One is an axiom which resembles the sufficiency axiom and indeed implies the latter. The second part specifies in more positive terms a mode of interpretation of statistical evidence, and thus resembles the confidence concept in some respects. \hfill ---\citep[p.\ 125]{birnbaum:1969}

}

The axiom that Birnbaum refers to says that two proportional likelihood functions contain the same evidence, but it is often written as a statement that the likelihood function contains the all of the evidence from data concerning the merits of various hypotheses. Each form of the axiom logically entails the other and so the difference of expression need not concern us. The second part that Birnbaum refers to is usually called the law of likelihood, and it says that the degree to which the evidence favors one hypothesis over another is given by the ratio of the likelihoods of the hypotheses. It is noted that in many accounts the likelihood principle is taken to be only what Birnbaum calls the axiom, and the law of likelihood is dealt with separately. However, for the purposes of the current paper it is necessary to follow Birnbaum's practice of including both parts in order that the critique of his alleged counter-example is fair.

 Birnbaum's own definitions of the likelihood principle are included verbatim in an appendix to this paper, but one is quite long and both are more technical than seems necessary. Their key features are captured within this succinct statement from Edwards:

\quote{
Within the framework of a statistical model, \textit{all} of the information which the data provide concerning the relative merits of two hypotheses is contained in the likelihood ratio of those hypotheses on the data, and the likelihood ratio is to be interpreted as the degree to which the data support one hypothesis against the other.  \hfill ---\cite[p. 31]{Edwards:1972:book}

}

That statement has a couple of shortcomings as a definition, aside from its lack of technical precision. First, the `hypotheses' dealt with in the likelihood principle are nothing more than values of the parameters of the statistical model, as Edwards demonstrates in examples throughout his book, and as Birnbaum makes clear in both of his definitions (see Appendix). That might not be thought to matter, as the use of  `hypothesis' as a synonym for `parameter value' is common practice in the discussion of likelihood functions and the likelihood principle---I've done so above---but there is scope for serious confusion. Not every hypothesis is a simple hypothesis corresponding to the value of a model parameter and not every possible hypothesis can be accommodated by any particular statistical model. Thus hypotheses will exist that cannot be compared using the likelihood principle. It is also important to note that the likelihood principle does not provide a mechanism for comparison of simple hypotheses that map onto a single value of a parameter with composite hypotheses that specify ranges of parameter values. Birnbaum writes explicitly  that the likelihood principle ``specifies no further structure or interpretation for the likelihood ratio scale, nor any specific concept of  `evidence supporting a set of parameter points.'\ '' \citep[p. 126]{birnbaum:1969}. Finally, the scope of the law of likelihood is naturally restricted because the scaling of a likelihood function is always arbitrary and entails an unknown proportionality constant. That unknown constant doesn't matter if two likelihoods of interest lie on the same likelihood function because the constant is the same for both and therefore cancels in the likelihood ratio, but otherwise it leads to potentially meaningless likelihood ratios. Thus the law of likelihood only tells us how to measure the relative support of parameter values that lie \textit{on the same} likelihood function. That restricted scope of the law of likelihood is not always noted, but Edwards includes it obliquely in his definition by writing ``Within a statistical model'', and it is implied in Birnbaum's definition where he writes that the law of likelihood ``consists of the statement that the evidence supporting one parameter point $i$ against another $i'$ is represented just by the numerical value of the likelihood ratio $L(i,i') = p_{ij}/p_{i'j}$'' because $i$ and $i'$ are two parameter values on the same likelihood function: ``$p_{ij}$ of $i, i \in \Omega$'' \citep[p.\ 126]{Birnbaum:1962} .

With all of that as a preamble, we can briefly state the likelihood principle in a manner that respects both Edwards's and Birnbaum's versions:

\bigskip
\noindent \textit{The likelihood principle} \; Two likelihood functions which are proportional to each other over the whole of their parameter spaces have the same evidential content. Thus, within the framework of a statistical model, all of the information provided by observed data concerning the relative merits of the possible values of a parameter of interest is contained in the likelihood function of that parameter based on the observed data. The degree to which the data support one parameter value relative to another on the same likelihood function is given by the ratio of the likelihoods of those parameter values.

\subsubsection{Example to illustrate the likelihood principle}

Consider a simple but artificial example where the probability of rain during each day of the current week are available via a seven day forecast from the local bureau of meteorology. Assume that there is no intermediate state between rain and not rain so that those states are mutually exclusive and exhaustive, and assume that the only time unit of interest is `day'. With those assumptions we get a statistical model where the probability distribution for each day consists of the forecast probability of rain and the complement of that as the probability of not rain. Using that statistical model we have only two possible observations, $x \in \{\text{rain},\text{not rain}\}$, and days of the week are the only model parameter, so $\theta \in \{\text{Monday, Tuesday, Wednesday, Thursday, Friday, Saturday, Sunday}\}$. Table \ref{tableRain} shows the probabilities. Each column of probabilities displayed in that table is a probability distribution that sums to unity, and the rows in the table are proportional to the likelihood functions associated with the observations of rain and not rain.

\bigskip
\begin{minipage}{\linewidth}
\centering
\captionof{table}{Probabilities of rain for the current week}
\label{tableRain}
  \begin{tabular}{rccccccc}
\toprule
 & Monday & Tuesday & Wednesday & Thursday & Friday & Saturday & Sunday \\
$P($rain$)$           & 0    & 0.07    & 0.65      & 0.2     & 0.05    & 0.01     & 0.01   \\
$P($not rain$)$       & 1    & 0.93    & 0.35      & 0.8     & 0.95    & 0.99     & 0.99  \\
\bottomrule    
\end{tabular}
\end{minipage}
\bigskip

There are two types of question that might be usefully addressed using the probabilities and likelihoods in the table. First, should you pack an umbrella? The answer will depend on your personal loss function (\textit{e.g.} do you care about a bit of rain?), but it might be reasonable to pack an umbrella on Wednesday and, maybe, Thursday. The answer to that umbrella question is informed by individual probabilities: the 65\% chance of rain on Wednesday provides a reason to carry an umbrella without reference to any other probability in the table. Such a decision treats the parameter value (Wednesday) as a given and the observable (rain) as a variable. The second type of question that can be addressed using the table goes the other way, from the observation to the parameter: what does an observation of rain tell you about what day it is? Say you wake from unconsciousness and observe that it is raining and wonder what day it is. You notice that table \ref{tableRain} is posted on the wall with a prominent label ``This week's rain probabilities''. The observation of rain means that only the upper row in table \ref{tableRain} is relevant, and within that row it is Wednesday that has the highest probability of rain occurring. The observation of $x=$ rain provides evidence in favor of the day being any day where there is a non-zero probability of rain, and it provides stronger evidence for the day being a day where there is a high probability of rain. Thus the observation $x=$ rain supports $\theta=$ Wednesday more strongly than it supports any other value of $\theta$. Wednesday is the maximal likelihood estimate of $\theta$ given the observation of rain and the probability model represented by table \ref{tableRain}.

It is worth stating explicitly that while the first type of question---should you carry an umbrella?---is addressed using individual probabilities, the second type of question---what day is it?---requires a comparison of all of the likelihoods in the relevant row of the table. The relevance of that distinction can be seen by considering the effect of changing the probability of rain in the table. If, for example, the probability of rain on Friday was changed from 0.05 to 1, it would not affect whether you should carry the umbrella on Wednesday, but it would change the maximal likelihood estimate from Wednesday to Friday even though the evidence itself, $x=$ rain, was not changed. The use of likelihoods to make inferences about parameter values on the basis of observations requires a \textit{comparison} of the likelihoods of the various values of the parameters.

The likelihood principle tells us that, given the model where probabilities of rain come from the forecast of the bureau of meteorology, all of the information from the observation of $x=$ rain relevant to estimation of the day of the week is contained in the likelihood function that is proportional to the probabilities in the top row of table \ref{tableRain}, and that the evidential favoring of that observation for each day relative to any other is proportional to the ratio of those likelihoods. Thus the observation of rain stands as evidence in favor of the day being Wednesday over the day being Thursday by the ratio of $0.65/0.2=3.25$ and simultaneously stands as evidence favoring Thursday over Saturday by the ratio $0.2/0.01=20$.
\section{Examples analogous to Birnbaum's problem}

Birnbaum's counter-example provides a situation where, at least on the face of it, application of the law of likelihood yields strong support for a hypothesis where intuitively it seems that there should be none.  The counter-example involves a single observation taken from a population of unknown mean and the hypotheses in question concern the spread of values in the population sampled: either it has a wide spread such that observations can take values of up to 100 either side of the unknown mean; or it has no spread so that all observations would take the value equal to the unknown mean. Birnbaum's analysis led him to say that, no matter what value is observed, the hypothesis that the population in question has no spread will have a markedly higher likelihood than the hypothesis that the population has a wide spread. 
 
Birnbaum's analysis is flawed, but the nature of the flaw is subtle and so it is helpful to begin with analyses of analogous examples that are more accessible than Birnbaum's in some details. Two examples are supplied here, one very simple and another slightly more complicated but directly analogous to Birnbaum's. Both utilize a model of the favored by statisticians: urns containing ed balls.

\subsection{Example 1: Known color}
\label{Ex1}

Imagine that an honest statistician places an urn on your desk and tells you that it is either urn 1a which contains 10{,}000 red balls or it is urn 1b which contains 100 red balls and 9{,}900 balls made up of 200 different colors, none of which is as numerous as the red balls.  He invites you to determine which urn you have on the basis of the colors of sampled balls. You may randomly sample as many balls as you like as long as it is sampling with replacement.

How many balls should you sample? The answer depends on what color(s) you observe in your sample. If you draw a ball of any  other than red then even a single ball is definitive evidence that the urn is 1b because only urn 1b contains non-red balls. If you draw a red ball then that is evidence in favor of the urn being urn 1a because 1a contains a larger proportion of red balls than urn 1b but, as urn 1b has \textit{some} red balls, the evidence of a single red ball would not be conclusive. Thus you might choose to draw more balls to increase the strength of evidence in favor of urn 1a over 1b. The larger the number of red balls observed without any other colors, the stronger that evidence would be, and more confident you could feel, that the urn was 1a, but drawing even a single non-red ball which would confirm that the urn was 1b. 

This is a simple problem that probably elicits the same intuitive response from all readers, one that matches that description in the previous paragraph. It is also one that matches the results of a formal likelihood analysis of the problem. 

A formal likelihood analysis of the problem requires a statistical model in which specifies the probabilities of the possible observations, and for this problem those probabilities obviously derive from the frequencies of the various colored balls within the urns. However we do not know most of those frequencies and so we cannot construct a complete sampling probability distribution based on individual colors. A convenient solution is available. As all non-red colors are equally informative about the question at issue (the identity of the urn), we can collapse the sample space into red and non-red to obtain a fully specified sampling distribution. At first glance, the relevant hypotheses concern the urn on your desk---$H_{1a}$: urn 1a; and $H_{1b}$: urn 1b---but the hypotheses in a likelihood analysis have to be parameter values. Conveniently we can obtain such hypotheses by noting that $H_{1a}$ implies that there is only one color of ball in the urn, and $H_{1b}$ there are 201 colors, so we have $H_{1a}': \;\nu_c=1$ and $H_{1b}': \; \nu_c=201$, where $\nu_c$ is a parameter specifying how many colors are in the urns. If $H_{1a}'$ is true then so is $H_{1a}$, and if $H_{1b}'$ is true then so is $H_{1b}$, so we can logically equate the hypotheses concerning the identity of the urn with values of the parameter $\nu_c$.

 Table \ref{tableEx1} shows the probabilities of observing red and non-red balls according to the relevant model for this problem. Note that the probabilities sum to unity across the rows because the two possible observations, $x=$ red and $x\neq$ red are exhaustive of the full sample space, but the columns do not sum to unity because the each row is a mutually exclusive condition or hypothesis.

\bigskip
\begin{minipage}{\linewidth}
\centering
\captionof{table}{Probabilities of drawing red and non-red balls for example 1}
\label{tableEx1}
\begin{tabular}{lcc}
\toprule
       & $x=$ red & $x \neq$ red \\ \cmidrule{2-3}
Urn 1a, \; $\nu_c=1$ & 1        & 0            \\
Urn 1b, \; $\nu_c=201$ & 0.01     & 0.99    \\
\bottomrule    
\end{tabular}
\end{minipage}
\bigskip

The relevant likelihood function is are derived from the probabilities in table \ref{tableEx1}. As the parameter values $\nu_c=1$ and $\nu_c=201$ are mutually exclusive and exhaustive, each likelihood function consists of only two points. For a sample of more than one ball we only need to keep track of the number of red and non-red balls because the sequence of the observations is not informative as a consequence of the sampling being with replacement. Thus, for a sample of $n=r+m$ balls consisting of $r$ red balls and $m$ non-red balls, the likelihood function is this:

\begin{equation}
\begin{aligned}
\label{LfnEx1alt}
L(H_{1a}')  &\propto P(x = \text{red} \,|\, H_{1a}')^{r} \times P(x \neq \text{red} \,|\, H_{1a}')^{m}\\
&\propto P(x = \text{red} \,|\, \nu_c=1)^{r} \times P(x \neq \text{red} \,|\, \nu_c=1)^{m} = 1^{r} \times 0^{m}\\
L(H_{1b}')  &\propto P(x = \text{red} \,|\, H_{1b}')^{r} \times P(x \neq \text{red} \,|\, H_{1b}')^{m}\\
&\propto P(x = \text{red} \,|\, \nu_c=201)^{r} \times P(x \neq \text{red} \,|\, \nu_c=201)^{m} = 0.01^{r} \times 0.99^{m}.
\end{aligned}
\end{equation}

\noindent $L(H_{1a}')=1$ where $m=0$ and is zero otherwise\footnote{Taking $0^{0}=1$.}. The likelihood ratio for any observed mixture of balls is just the ratio of those two values. 

\begin{equation}
\label{LratioEx1alt}
\frac{L(H_{1a}')}{L(H_{1b}')} = \frac{1^r \times 0^m}{ 0.01^{r} \times 0.99^{m}}
\end{equation}

Birnbaum's problem involves only a single observation, so it is worth calculating the likelihood ratio for this example under the two possible $n=1$ outcomes. If a red ball is drawn, $r=1$ and $m=0$, the likelihood ratio is 100:1 in favor of $H_{1a}'$, and if a non-red ball is drawn, $r=0$ and $m=1$, the likelihood ratio is 0, which is conclusive evidence in favor of $H_{1b}'$ over $H_{1a}'$. That likelihood analysis matches the intuitive responses to the possible results in this example, as described above.

\subsection{Example 2: Unknown color}
\label{Ex2}

The second example is identical to the first except that the color shared by the two urns is not known. Thus, imagine that an honest statistician places an urn on your desk and tells you that it is either urn 2a which contains 10{,}000 balls of one color, or it is urn 2b which contains 100 balls of that same color and 9{,}900 balls made up of 200 different colors, none of which is as numerous as the balls with the shared but unknown color.  He invites you to determine which urn you have on the basis of observed colors of sampled balls. You may sample as many balls as you like as long as it is sampling with replacement.

How many balls should you sample? As in example 1, the answer depends on what color(s) are observed. This is more complex than the previous example, but an intuitive interpretation of possible outcomes is still fairly clear. This time there is no possible $n=1$ observation that could stand as evidence for one urn over the other because any single color might be the color that is present in both urns, and a sample of one ball will certainly have a single color. Thus the sample needs more than one ball for it to provide evidence that discriminates between the urns. If all of the balls in a sample of more than one have the same color, that would stand as evidence in favor of the urn being 2a. However, as soon as the sample contains multiple colors the evidence conclusively favors urn 2b.

As in the previous example, the statistical model has probabilities determined by the frequencies of the ball colors in the urns, but in this example we do not know \textit{any} of those frequencies. The probabilities can only be had by making an assumption about the color that is shared by the urns, a color that will be denoted by $\mu$. For example, if we assume that $\mu=$ blue, the probabilities are those in table \ref{tableEx2}. That table is larger than the table in example 1 because the model for this example has two unknown parameters, $\mu$ and $\nu_c$, where the model in example 1 only has $\nu_c$. The left and right halves of table \ref{tableEx2} are mutually exclusive alternatives, so each row in the table contains two probability distributions and sums to 2. Unfortunately, the probabilities in table \ref{tableEx2} do not yield a usable likelihood function, even if an observed color is substituted for 'blue' in the table, as we would not know which half of the table is relevant.

\bigskip
\begin{minipage}{\linewidth}
\centering
\captionof {table}{Probabilities of drawing blue and non-blue balls for example 2}
\label{tableEx2}
\begin{tabular}{@{}llcclcc@{}}
\toprule
       &  & \multicolumn{2}{c}{$\mu=$ blue} &  & \multicolumn{2}{c}{$\mu \neq$ blue}                              \\ \cmidrule{3-4} \cmidrule{6-7}
       &  & $x=$ blue     & $x \neq$ blue    &  & \multicolumn{1}{l}{$x=$ blue} & \multicolumn{1}{l}{$x \neq$ blue} \\ \midrule
Urn 2a, \;$\nu_c=1$ &  & $1$          & $0$             &  & $0$                          & $1$                              \\
Urn 2b \;$\nu_c=201$ &  & $0.01$       & $0.99$          &  & $< 0.01$                  & $> 0.99$                      \\ \bottomrule
\end{tabular}
\end{minipage}
\bigskip

A likelihood function relevant to this problem can be obtained, but its derivation is less obvious than in the previous example. The difficulty is that even though we don't really care about the value of $\mu$ as it is not relevant to the question about the identity of the urn, without knowledge of $\mu$ we cannot interpret the observation of any single color. However, as the question of interest does not concern the color of balls, a likelihood function based on the \textit{number} of colors in the sample rather than the \textit{nature} of the color(s) can serve our purpose. Conveniently, that  function is independent of the nuisance parameter $\mu$. 

Let $n_c$ be the number of distinct colors in a sample of $n$ balls. Table 3 presents the probabilities for observations of $n_c=1$ and $n_c>1$. The value $(\leqslant 0.01)^{n-1}$ in that table requires a little explanation. The probability of observing $x=\mu$ when a ball is drawn from urn 2b is $100/10{,}000=0.01$, and because of the other colors are less numerous (for consistency with Birnbaum's example, below), the probability of observing any other color is less than 0.01. The exponent is $n-1$ because even thought any observation would satisfy $n_c=1$ when $n=1$, repeated draws of a single color are needed for $n_c=1$ with $n=2$ or more. 

\bigskip
\begin{minipage}{\linewidth}
\centering
\captionof{table}{Probabilities of a sample of $n$ balls containing $n_c$ distinct colors in example 2}
\label{tableEx2b}
\begin{tabular}{lcc}
\toprule
       & $n_c=1$ & $n_c>1$ \\ \cmidrule{2-3}
Urn 2a, \;$\nu_c=1$ & 1        & 0            \\
Urn 2b, \;$\nu_c=201$ & $(\leqslant 0.01)^{n-1}$     & $1- (\leqslant 0.01)^{n-1}$   \\
\bottomrule    
\end{tabular}
\end{minipage}
\bigskip

We use the same logic as in the previous example to equate the hypotheses about the identity of the urns, $H_{2a}$ and $H_{2b}$, with hypotheses $H_{2a}'$ and $H_{2b}'$ which relate to the parameter $\nu_c$. For an observed  $n_c=1$ from a sample of $n$ balls we get the likelihood function

\smallskip
\begin{equation}
\begin{aligned}
\label{LfnEx2}
L(H_{2a}') &\propto P(n_c=1 \,|\, \nu_c=1)  = 1\\
L(H_{2b}') &\propto P(n_c=1 \,|\, \nu_c=201)  =  (\leqslant 0.01)^{n-1},
\end{aligned}
\end{equation}

\noindent and for $n_c>1$ the likelihood function is
\smallskip
\begin{equation}
\begin{aligned}
\label{LfnEx2ncgtn}
L(H_{2a}') &\propto P(n_c > 1 \,|\, \nu_c=1)  =  0\\
L(H_{2b}') &\propto P(n_c > 1 \,|\, \nu_c=201)  =  1- (\leqslant 0.01)^{n-1}.
\end{aligned}
\end{equation}

\noindent The likelihood ratio also depends on whether the observation is $n_c=1$ or $n_c>1$. For $n_c=1$ it is
\smallskip
\begin{equation}
\label{ratioEx2n1}
\frac{L(H_{2a}')}{L(H_{2b}')} =\frac{1}{(\leqslant 0.01)^{n-1}},
\end{equation}

\noindent and for $n_c>1$ it is
\smallskip
\begin{equation}
\frac{L(H_{2a}')}{L(H_{2b}')} =\frac{0}{1-(\leqslant 0.01)^{n-1}}.
\end{equation}

Again we should consider the possible results when only a single ball is drawn for consistency with Birnbaum's problem. If $n=1$ then the test statistic will necessarily be $n_c=1$ so the relevant likelihood function is equation \ref{ratioEx2n1}. Where $n=1$,  $n-1 = 0$, so the likelihood ratio will be 1:1. That result means that no matter what color is observed with a single draw the evidence is neutral with respect to the hypotheses concerning the identity of the urn. That result matches an intuitive response to the example.

\section{Birnbaum's counter-example}

Birnbaum's example is similar to example 2 above, but in the space of numbers rather than colors. It consists of a single observation, $x$, drawn from a population with unknown mean $\mu$ and a spread parameter, $\sigma$, that is either 0 or 100. Birnbaum writes\footnote{There are several obvious typographical errors in the original that have been corrected here to avoid confusion. None of the errors changed the sense of the text or the formulae.}:

\begin{quote}
Let $\theta=(\mu,\sigma)$: The values $\sigma=0$ or $100$, respectively, represent the unknown precision, either very high or very low, of a single observation (measurement) $x$; and $\mu$, an integer, $-10^{10}<\mu<10^{10}$, is the unknown true value of a quantity to be measured. Let 

\end{quote}

\begin{small}
\begin{equation*}
f(x,\theta) = f(x,\mu,\sigma)=
\begin{cases} 1, \text{ if } \sigma=0 \text{ and }x=\mu\\
0, \text{ if } \sigma=0 \text{ and } x\neq \mu, \text{ for all }x,\mu; \\
c[100-|x-\mu|] \text{ if } \sigma = 100, \text{ for }|x-\mu|<100,\\
0 \text{ otherwise.}
\end{cases}
\end{equation*}

\end{small}

\begin{quote}
 Thus $\sigma=0$ gives an error-free measurement $x=\mu$; $\sigma = 100$ gives a wide ``triangular'' distribution of $x$ centered at $\mu$. The constant $c=1/10{,}040 \doteq 10^{-4}$ (to give total probability 1). The likelihood function determined by any of the possible observations $x, |x| < 10^8$, is
 
\end{quote}
\begin{small}
\begin{equation*}
\begin{split}
f(x,\mu,\sigma)
& =1, \text{ for } \sigma=0 \text{ and } \mu=x; \\
  & = 100c \doteq 10^{-2}, \text{ for } \sigma=100 \text{ and } \mu = x;\\
& <100c \doteq 10^{-2}, \text{ for } \sigma = 100 \text{ and } \mu \neq x.
\end{split}
\end{equation*}
\end{small}

\begin{quote}
\noindent A simple application  of the likelihood principle, in accord with the account of the preceding section\footnote{[Reproduced in the appendix to this paper.]}, is the following: The parameter point $\theta = (\mu,\sigma) = (x,0)$, that is, $\mu=x$ and $\sigma=0$, has uniquely maximum likelihood. (Its likelihood ratio, as against each other point in turn, exceeds unity.) Thus it seems to be supported by the observed result, as against each and all other parameter points. This includes, in particular, that the value $\sigma=0$ seems to be supported against the alternative $\sigma=100$. Evidently the numerical values of the likelihood ratios referred to, all exceeding 100, represent strong evidence in some sense. Evidently such interpretations of statistical evidence must be regarded as misleading, and strongly misleading, in case the true value of $\sigma$ is not 0 but 100. But such misleading interpretations will be suggested by the likelihood principle with probability unity, if for example $\sigma=100$ and $\mu=0$ are the true parameter values, since in this case each possible outcome $x$ determines a likelihood function of the form considered. \hfill ---\citep[pp. 127--128]{birnbaum:1969}

\end{quote}

Birnbaum's analysis leads to the conclusion that no matter what single value of $x$ is observed, the likelihood associated with the hypothesis that $\sigma = 0$ will always be at least 100 times larger than that associated with the hypothesis $\sigma=100$, even if the latter hypothesis is true. It is generally assumed that the example is therefore a counter-example to the likelihood principle.

\section{Critique}


Birnbaum's alleged counter-example seems to be closely analogous to example 2 presented in section \ref{Ex2} in that they both present a choice between two populations which have a fixed, but unknown, parameter in common and a parameter representing distributional spread that takes one of two specified values. However, despite the structural similarities, example 2 and the alleged counter-example yielded quite different outcomes---one seems to support the likelihood principle and one that is generally taken to represent an important counter-example to the likelihood principle. At least one of three things must be true: the examples are disanalogous in some important way; one of the analyses is flawed; the counter-example is a counter-example to the particular analysis used by Birnbaum rather than to the likelihood principle in general. Those possibilities are be discussed in sections below.

\subsection{The examples are disanalogous}

There are several reasons to conclude that there is a close analogy between example 2 and Birnbaum's example. 

\begin{enumerate}

\item The desired inferences in the examples are equivalent in that both offer a choice between two hypothetical populations. 

\item One population in each example consists of a single class of objects whereas the other has more than one.  

\item Both examples use discrete scale of measurement so that the colors of the balls can be mapped one to one onto integer values in Birnbaum's example.

\item The parameters of the statistical models can be mapped one to one onto each other: $\mu$ is a value shared by the two possible population in both; and Birnbaum's spread parameter, $\sigma$, has its equivalent in example 2, the number of colors of ball in the urns, $\nu_c$, with the simple linear relationship $(\nu_c-1)/2 = \sigma$.  Thus the parameter space for both examples is equivalent and can be denoted in the same manner: $\theta=(\mu, \sigma)$.

\item \label{twoObs}In both Birnbaum's example and example 2 a second observation would yield either definitive evidence that $\sigma=100$ if the two observations were dissimilar, or strong evidence with a likelihood ratio of at least $10^4$ in favor of $\sigma=0$ if the observations were identical.


\end{enumerate}

Those considerations seem sufficient to conclude that example 2 is a close analogue of Birnbaum's example.


\subsection{One of the analyses is flawed}
\label{flawed}

If Birnbaum's example and example 2 are closely analogous then the divergent results must come from differences in the analyses used. That can be confirmed by simply applying the analysis used in example 2 to Birnbaum's problem. Therefore, assume that we observe $x=17$ so that the full set of probabilities are those in table \ref{tableBBProbs}. That table can be obtained by just substituting  into the example 2 table `17' for `blue' and approximating $100c$ with $0.01$, and, for a likelihood analysis, it presents exactly the same quandary as it did in example 2: which half of the table is relevant? Using the number of distinct values in the sample, $n_v$, as a test statistic side-steps that difficulty and gives the probabilities in table \ref{tableBBexShort}. Those probabilities are the same as those in table \ref{tableEx2b}, so the likelihood functions are also the same as those in example 2. That demonstrates that application of the example 2 analysis to Birnbaum's problem yields a result where the likelihood-based interpretation of the evidential meaning of the observation matches the intuitive response. It can be concluded that the essence of the counter-example lies in the analysis rather than the logical structure of the problem.

\bigskip
\begin{minipage}{\linewidth}
\centering
\captionof{table}{Probabilities of observing $x=17$ and $x \neq 17$ for Birnbaum's example (using $100c = 0.01$).}
\label{tableBBProbs}
\begin{tabular}{@{}llcclcc@{}}
\toprule
            &  & \multicolumn{2}{c}{$\mu=17$}              &  & \multicolumn{2}{c}{$\mu \neq 17$}                            \\ \cmidrule{3-4} \cmidrule{6-7}
            &  & $x=17$             & $x \neq 17$          &  & \multicolumn{1}{l}{$x=17$} & \multicolumn{1}{l}{$x \neq 17$} \\ \midrule
$\sigma=0$  &  & $1$                & $0$                  &  & $0$                        & $1$                             \\
$\sigma=100$ &  & $0.01$ & $0.99$ &  & $<0.01$                    & $> 0.99$                     \\ \bottomrule
\end{tabular}
\end{minipage}
\bigskip

\bigskip
\begin{minipage}{\linewidth}
\centering
\captionof{table}{Probabilities of a sample of $n$ observations containing $n_v$ distinct values in Birnbaum's example (using $100c = 0.01$).}
\label{tableBBexShort}
\begin{tabular}{lcc}
\toprule
       & $n_v=1$ & $n_v>1$ \\ \cmidrule{2-3}
$\sigma = 0$ & 1        & 0            \\
$\sigma = 100$ & $(\leqslant 0.01)^{n-1}$     & $ 1- (\leqslant 0.01)^{n-1}$   \\
\bottomrule    
\end{tabular}
\end{minipage}
\bigskip


%

\subsection{Objections to Birnbaum's analysis}

The success of the example 2 analytical formulation in dealing with Birnbaum's example does not necessarily mean that Birnbaum's analytical formulation is flawed, but 
 two interrelated objections can be raised. First, Birnbaum's analysis leads to the inference being determined by an inferentially irrelevant parameter and, second, the statistical model used in that analysis has two parameters, $\mu$ and $\sigma$, whereas the problem supplies only a single observation.

\subsubsection{Objection 1: outcome determined by inferentially irrelevant parameter}

Birnbaum's presentation of the parameter of interest is somewhat contradictory. In the description of the problem he writes that $\mu$ is ``the unknown true value of a quantity to be measured'', but that only means that $\mu$ is the parameter of interest to the agent making the measurement. That agent does not face any monstrous dilemma in determining the best supported value of that parameter because there is nothing misleading in the support of $\mu=x$ over alternatives that the likelihoods provide. The allegedly misleading evidence concerns not $\mu$, but $\sigma$, as Birnbaum makes clear:

\quote{
This includes, in particular, that the value $\sigma=0$ seems to be supported against the alternative $\sigma=100$. [\dots] Evidently such interpretations of statistical evidence must be regarded as misleading, and strongly misleading, in case the true value of $\sigma$ is not 0 but 100. But such misleading interpretations will be suggested by the likelihood principle with probability unity, if for example $\sigma=100$ and $\mu=0$ are the true parameter values, since in this case each possible outcome $x$ determines a likelihood function of the form considered.
}

As the alleged counter-example rests entirely on the behavior of the analysis with regard to $\sigma$, the parameter of interest to us must be $\sigma$. In that case the value of $\mu$ is incidental to the problem---in so far as it plays a role in the analysis it is as a nuisance parameter. An appropriately designed analysis would ensure that estimation of $\sigma$ is not influenced by that nuisance parameter. Section \ref{flawed} contains such an analysis, and it provides a likelihood function supporting conclusions that mirror our intuitive response to evidence, in stark contrast to those that are drawn from Birnbaum's function.

We have two different likelihood functions that are, apparently, appropriate for the same experiment, and they yield contrasting pictures of the evidence. The question therefore arises as to which of those functions is correct. They are both correct in the sense of being valid likelihood functions for the experiment, despite Birnbaum's being incomplete\footnote{Birnbaum's likelihood function lacks a likelihood for $\theta=(\mu\neq x,\sigma=0)$. As the likelihood of that parameter value is zero, its omission does not much matter.}. However, each function has a different inferential scope. Birnbaum presents the likelihood function for the vector parameter $\theta=(\mu,\sigma)$ and it provides likelihood ratios that quantify the relative merits of the various values of $\theta$, in particular it shows that any observation of $x$ supports $\theta=(x,0)$ over any other value of $\theta$. It may seem natural to assume that the best supported value of $\sigma$ is therefore 0, but that assumption would be wrong. Simultaneous assessment of the evidence regarding $\mu$ and $\sigma$ as components of a vector need not yield the same result as the assessment of the evidence regarding either of them as independent scalars. In the context of Birnbaum's problem, the $\mu$ dimension of $\theta$ affects the probability of observing $x$ more sharply than does the $\sigma$ dimension and so the relationship between $x$ and $\mu$ plays a dominant role in shaping the evidence regarding values of $\theta$. Of course, we are interested in $\sigma$ rather than $\mu$ or $\theta$, so a dominant influence of $\mu$ over $\sigma$ in the analysis is disastrous. We need to deal with $\mu$ as a nuisance parameter. In sections \ref{Ex2} and \ref{flawed} $\mu$ is eliminated by the choice of parameterization, but there are other strategies could be equally effective, in particular a marginal likelihood function for $\sigma$ obtained by integrating out $\mu$ would yield the same result of equal support for $\sigma=0$ and $\sigma=100$. 
 Thus there is no contradiction in the data providing equal support for $\sigma=0$ and $\sigma=100$ at the same time as it shows strong support for $\theta=(x,0)$ over $\theta=(x,100)$, but where the parameter of interest is $\sigma$ the latter result is not relevant.


If Birnbaum's example is a counter-example, it is a weak counter-example in that the misleading outcome is very easily avoided by appropriate treatment of a nuisance parameter.

\subsubsection{Objection 2: more parameters than data}
\label{overfitting}

Birnbaum's model has more parameters than the number of available data points, or, given the vector notation used in places for the parameter $\theta$, the model has a parameter with more dimensions than the number of data points. Either way it represents a condition that is often called `overfitting'. The pivotal role of overfitting in Birnbaum's example is readily seen by simply considering two variants of his example that lack overfitting: (i) where there is only one parameter; and (ii) where there is a second datum. 

\begin{description}
\item{\textit{(i) Single parameter variant}} \quad Say that it is known that the true value to be measured is $\mu=17$ and we wish to determine whether $\sigma=0$ or $\sigma=100$ in order to know about the measurement precision. The model has only one parameter to be estimated from the data and the relevant probabilities are those in the left half of table \ref{tableBBProbs}. This single parameter variant of Birnbaum's problem has exactly the structure of the balls problem presented in example 1, and the intuitive response to the evidence presented by any possible observation matches exactly the relevant likelihood function. Thus the alleged counter-example does not survive modification into a single parameter variant.

\item{\textit{(ii) Two observation variant}}\quad Say that neither part of $\theta=(\mu,\sigma)$ is known and two values are observed, $x_1$ and $x_2$. It should be intuitively clear that if $x_1=x_2$ then there would be strong support for the hypotheses that $\theta=(\mu=x_1=x_2,\sigma=0)$ and, equivalently, if $x_1$ and $x_2$ differ then they would provide proof that $\sigma\neq0$ and thus support $\theta=(\frac{x_1+x_2}{2},\sigma=100)$ over any other value of $\theta$. The type of analysis used in example 2, where the number of distinct observations, $n_c$, is used as a test statistic, would be applicable to this problem, but an analysis using the observed values directly is equally applicable and, as it is more similar to the original analysis by Birnbaum, it will be presented, first for the case where $x_1=x_2=17$ and then where $x_1=17, x_2=132$.

For an observation of $x_1=x_2=17$, the likelihood function consists of two lines on $\mu$:

\begin{align}
 L(\mu,\sigma=0) &\propto 
 \begin{cases} 1, \;\text{where } |17-\mu|=0\\
0, \; \text{ otherwise}
\end{cases} \\
L(\mu,\sigma=100) &\propto 
\begin{cases} (c(100-|17-\mu|))^2, \;\text{where } |17-\mu| < 100 \\ 0, \;\text{ otherwise}.
 \end{cases}	
\end{align}

That likelihood function supports the parameter $\theta=(17,0)$ over all others by a margin of at least $10^4$, but, given that an observation of $x_1=x_2$ would occur with a frequency of less than $10^{-4}$ if $\sigma = 100$, that support would very rarely be misleading.

For an observation of $x_1=17, x_2=132$ the likelihood function is also two lines on $\mu$, but one of the lines, where $\sigma =0$, is a flat line at zero and the other, where $\sigma=100$, consists of a line mostly at zero but with a central concave down parabola in the region between $\mu=33$ and  $\mu=117$:

\begin{align}
 L(\mu,\sigma=0) &\propto \begin{cases} 1, \;\text{where } |17-\mu|=0\\
0, \; \text{ otherwise}
\end{cases} \\
L(\mu,\sigma=100) &\propto \begin{cases} (c(100-|17-\mu|))^2, \;\text{where } |17-\mu| < 100 \\ 0, \;\text{ otherwise}. 
\end{cases}
\end{align}

That likelihood function displays strongest support for the parameters $\theta=((x_1+x_2)/2, 100)$, and that support would never be misleading concerning the $\sigma$ dimension of $\theta$.

Thus it is demonstrated that the likelihood analysis of the problem when there are two data points results in strong support for $\theta$ containing $\sigma=0$ when the two observations are identical and $\theta$ containing $\sigma=100$ when they differ, both cases which are the same as the intuitive response to the problem. The probability of obtaining strong misleading support for $\sigma=0$ when $\sigma$ is actually 100 is less than $10^{-4}$. The alleged counter-example does not survive the inclusion of a second observation.

\end{description}

Both modifications of Birnbaum's example to allow the number of data points to equal the number of estimated parameters disarmed the counter-example by making probability of misleading support for a false hypothesis very low or zero. Thus, to the extent that Birnbaum's example is a counter-example, it is a counter-example to the likelihood principle only in the specific circumstance that the number of parameters exceeds the number of data points. Such reliance on overfitting may of itself be sufficient to dismiss the alleged counter-example because many statistical approaches are known to behave erratically where the number of estimated parameters exceeds the number of data points, and such misbehavior is not routinely assumed to be indicative of a counter-example to a statistical or philosophical principle.

\section{What happened had to happen}

The assumption that the likelihood principle allows the best supported hypothesis in any investigation to be a determinist hypothesis---a hypothesis that whatever happened had to happen---stands as an important criticism of the likelihood principle. Two distinct types of determinist hypotheses occur. The first kind corresponds to a (set of) parameter value(s) within a statistical model that provide for only a single specific outcome. The second kind of determinist hypothesis is that no matter what happened, it had to happen. Each kind is discussed in this section.

 A determinist hypothesis of the first kind is instantiated in Birnbaum's analysis by the parameter value $\theta=(x,0)$: no matter what single value of $x$ is observed,  that parameter value has a higher likelihood than any other even if the true value is $\theta=(y,100)$. It leads Birnbaum to write that ``misleading interpretations will be suggested by the likelihood principle with probability unity'', but that is a true statement in only a very limited sense. If it is intended to mean that misleading interpretations will be suggested in \textit{every} likelihood analysis of  \textit{every} investigation, it is false. Neither the non-overfitted variants of Birnbaum's problem set out in section \ref{overfitting} nor the urn problems presented in sections \ref{Ex1} and \ref{Ex2} contain parameter values that would be misleadingly supported with unity probability. Thus the statement that ``misleading interpretations will be suggested by the likelihood principle with probability unity'' is relevant specifically to the analysis and experiment that Birnbaum supplies. The presence of this type of determinist hypothesis in a statistical model is not universal and may be rare.
 
 The universality of the whatever happened had to happen hypothesis comes only with the second type of determinist hypothesis, one that Jaynes described as a `sure thing' hypotheses wherein ``every detail of the data was inevitable'' \citep[p.\ 195]{Jaynes:2003}. To obtain such a hypothesis requires invocation of Descartes' evil demon, a device beloved by philosophers. The evil demon is a creature who is able to determine the value of observations and is determined to mislead the investigator, and \citet[pp. 183--187]{Sober:1988} argues that the evil demon is a nuisance parameter in every likelihood analysis even though it is usually invoked only implicitly. For example, \citet{mayo_and_spanos:2011} criticise the likelihood principle using a very straightforward example where there is a determinist hypothesis, $H_i^*$, in the statistical model for a coin-tossing experiment.
 
\quote{
For an extreme case, $H_i^*$ can assert that the probability of heads is 1 just on those cases that yield heads, 0 otherwise. [\dots] So the fair coin hypothesis is always rejected in favor of $H_i^*$, even when the coin is fair. \hfill---\citep[p.\ 184]{mayo_and_spanos:2011}

}

 There is no mention an evil demon, but it would be an evil demon who sets the probability of heads for each toss of the coin in $H_i^*$. That evil demon provides a reason to mistrust the likelihood principle because everyone would agree that a preference for $H_i^*$ on the basis of \textit{any} observed data would be silly. However, contrary to what might be implied by Mayo \& Spanos, the likelihood principle does not entail such a preference because the fair coin hypothesis and $H_i^*$ lie on distinct likelihood functions because they exist in different parameter spaces. Not only that, but they cannot be contained in a single simple statistical model because of their probability models are contradictory. To flesh that out, not that the statistical model universally used for a coin tossing experiment has a parameter space consisting of $p$, the fixed but unknown probability of heads for each trial, and $n$, the number of trials. The probability of any observable number of heads, $h$ where $0 \leqslant h \leqslant n$, is obtained from the binomial distribution, $Bin(n,p)$, and so the likelihood function for $p$ when we observe $h$ heads is:

\begin{equation}
\label{binom}
L(p) \propto p^h (1-p)^{n-h}
\end{equation}

\noindent That model treats the value of $p$ is an unknown constant that is the same for each toss of the coin, and the fair coin hypothesis maps onto the hypothesis that $p=1/2$. Contrast that with Mayo \& Spanos's clearly stated determinist hypothesis whereby $p$ varies from toss to toss. The varying $p \in \mathbb{N}\{0,1\}$ of the determinist hypothesis is a different thing from the fixed but unknown constant $p \in \mathbb{R}, 0 \leqslant p \leqslant 1$, of the set of hypotheses that include the fair coin hypothesis, and so it is not among the points of the likelihood function given in equation \ref{binom}. As the likelihood principle offers no guidance as to how to compare the nonsensical $H_i^*$ with the fair coin hypothesis, the likelihood principle cannot be blamed for any misguided acceptance of $H_i^*$. Further, in this context, a lack of guidance should hardly be seen as a failing.

Some will be tempted to say that the determinist hypothesis can be grafted onto the likelihood function given in equation \ref{binom} and therefore that hypothesis can serve to undermine the likelihood principle even where it is not explicitly included in the function. However, such a cobbled together model is inappropriate for two reasons. First, the both the parameter spaces and the probability mechanisms of the two parts of the statistical model would be distinct and so a cobbled together union of two models remains two distinct models, one where the outcomes result from probabilities determined by the binomial distribution using the parameters $n$ and $p$, and the other where the outcomes are determined by the whim of an evil demon. In that case, some prior probabilities would be necessary to weight the two parts of the combined model and, as others have pointed out \citep{Edwards:1992, pawitan:2001:book, royall:2004}, one would inevitably lean towards a very low prior probability on the evil demon arm of the model. Second, given that the choice of statistical model is a arbitrary decision on the part of the analyst, why would anyone chose to include an evil demon hypothesis along with the more interesting hypotheses? It is no more sensible to consider the evil demon hypothesis than it would be to consider, for example, a probability model for the urn problems where the probability assigned to the drawing of green balls is 7.9 times the proportion of balls in the urn that are green. Even though it could be done, it would not be sensible to do so.

The whatever happened had to happen hypothesis has only a very limited capacity to undermine the likelihood principle. It is relevant to likelihood analyses where there is overfitting, as in Birnbaum's analysis but not otherwise, and, even so, it is probably no more important to avoid overfitting in likelihood analyses than it is in other statistical analyses. The evil demon nuisance parameter is easily avoided by simply choosing a statistical model which does not include it, as all sensible models probably do. No example where an evil demon is invoked, explicitly or implicitly, can serve as a general counter-example to the likelihood principle. The notion that a likelihood analysis will \textit{always} support a hypothesis that whatever happened had to happen is wrong.

\section{Appendix: The likelihood principle according to Birnbaum}
This appendix is the contains two versions of the likelihood principle given by Birnbaum, the first from his 1962 paper which contains his proof that the likelihood principle is entailed by the conditionality and sufficiency principles \citep{Birnbaum:1962} and the second from his paper of 1969 which contains the alleged counter-example that is the subject of this paper \citep[pp. 125-126]{birnbaum:1969}.

\subsection{Birnbaum's 1962 version}

\begin{quotation}

\begin{small}
	\noindent \textit{The likelihood principle (L):} If $E$ and $E'$ are any two experiments with the same parameter space, represented respectively by density functions $f(x,\theta)$ and $g(y,\theta)$; and if $x$ and $y$ are any respective outcomes determining the same likelihood function; then Ev$(E,x)$ = Ev$(E',y)$. That is, the evidential meaning of any outcome $x$ of any experiment $E$ is fully characterized by giving the likelihood function $cf(x,\theta)$ (which need be described only up to an arbitrary positive constant factor), without reference to the structure of $E$.
	\end{small}
	\end{quotation}

\subsection{Birnbaum's 1969 version}

\begin{quotation}
\begin{small}
Another general concept of statistical evidence, the \textit{likelihood concept} (often called the \textit{likelihood principle}), consists of two parts. One is an axiom which resembles the sufficiency axiom and indeed implies the latter. The second part specifies in more positive terms a mode of interpretation of statistical evidence, and thus resembles the confidence concept in some respects.

To formulate these, we require the definition of the \textit{likelihood function} (which plays important roles, technical and conceptual, in various areas of mathematical statistics): For each model $(E,j)$ of statistical evidence, the function $p_{ij}$ of $i, i \in \Omega$, is called the likelihood function. Here $j$ is fixed. More precisely, the function $p_{ij}$ is specified as one among many alternative, equivalent representations of the same likelihood function, all having the form $cp_{ij}$ where $c$ denotes an arbitrary positive number. 

We discuss the axiom first:

\begin{description}
\item[(L):]
The likelihood axiom: If two models of statistical evidence ($E,j$) and ($E',j'$) determine the same likelihood function, then they represent the same evidential meaning. That is, if for some positive $c$ we have $p_{ij}=cp'_{ij'}$ for each $ i \in \Omega$ (the common parameter space of $E$ and $E'$), then $Ev(E,j) = Ev(E',j')$ (where $E = (p_{ij}), E' = (p'_{ij'})$).
\end{description}

We note that when we take $E' = E$, this becomes just the sufficiency axiom.

The second part of the likelihood concept complements the likelihood axiom, by indicating to some extent how a likelihood function may be interpreted as a representation of evidential meaning. It consists of the statement that the evidence supporting one parameter point $i$ against another $i'$ is represented just by the numerical value of the likelihood ratio $L(i,i') = p_{ij}/p_{i'j}$, with the value  of unity marking neutral evidence and successively larger values indicating stronger support for $i$ against $i'$. The concept specifies no further structure or interpretation for the likelihood ratio scale, nor any specific concept of  ``evidence supporting a set of parameter points.''

\end{small}
\end{quotation}


\bibliographystyle{agsm}
\bibliography{bib2}

\end{document}